\documentstyle{article}
\input epsf.tex

\title{Flexible polyhedra in the Minkowski 3-space}

\author{Victor Alexandrov}

\date{Sobolev Institute of Mathematics, Novosibirsk-90, 630090, Russia.
E-mail: alex@math.nsc.ru}

\begin{document}

\maketitle

\begin{abstract}
We prove that flexible polyhedra do exist in the Minkowski
3-space and each of them preserves the (generalized) volume and the (total) mean
curvature during a flex. To prove the latter result, we introduce the notion of the
angle between two arbitrary non-null nonzero vectors in the Minkowski plane.


\end{abstract}

\section{Introduction}

Recall that the Minkowski $n$-space ${\bf R}^n_1$ is the linear space
which consists of all ordered $n$-tuples of reals $x=(x_1,x_2,\dots ,x_n)$ and is endowed
with the following scalar product: $(x,y)=x_1y_1+x_2y_2+\dots +x_{n-1}y_{n-1}-x_ny_n$.
The length $|x|=\sqrt{(x,x)}$ of a vector $x$ is either a positive number or
the product of a positive number by the imaginary unity $i$, or zero
(see, for example, \cite{14}).

Let $\Sigma$ be a connected $(n-1)$-dimensional simplicial complex which is a
manifold. A continuous map $P:\Sigma\to {\bf R}^n_1$ is said to be a polyhedron in ${\bf R}^n_1$ if
it is affine and injective on each simplex of $\Sigma$. However we call the image
$P(\Sigma)$ a polyhedron too. We say that a two-dimensional  polyhedron $P$ is disk- or
sphere-type if the body of $\Sigma$ is homeomorphic to a disk in
the Euclidean 2-space or a round sphere in the Euclidean 3-space, respectively.

A polyhedron $P=P(\Sigma)$ is said to be flexible if there exists an analytic
(with respect to a parameter) family of polyhedra $P_t=P_t(\Sigma)$ $(0\leq
t\leq 1)$ such that

i) $P_0=P$;

ii) for each edge $e$ of $\Sigma$, the length of $P_t(e)$ is independent of $t$;

iii) there exist two vertices $v_1$ and $v_2$ of $\Sigma$ such that the
(Minkowski) distance between $P_t(v_1)$ and $P_t(v_2)$ is nonconstant in $t$.

The family $P_t$ $(0\leq t \leq 1)$ is called a nontrivial flex of $P$.

Similarly, we can define the notion of a flexible polyhedron in a Euclidean (as
well as spherical and hyperbolic) space.

The most important results in the theory of flexible polyhedra in the Euclidean
3-space are as follows:

a) there are no convex flexible polyhedra (A.Cauchy, 1813 \cite{6};
A.D.Alek\-sand\-rov, 1950 \cite{1});

b) there exist flexible embedded sphere-type polyhedra (R.Connelly, 1977
\cite{7}, \cite{12}; K.Steffen, 1980 (see, for example, \cite{3}));

c) each flexible polyhedron preserves its mean curvature during a flex
(R.Alexander, 1985 \cite{2}; F.Almgren and I.Rivin, 1998 \cite{4});

d) each flexible polyhedron preserves its (generalized) volume during a flex
(I.Kh.Sabitov, 1996 \cite{15}--\cite{17}; R.Connelly, I.Sabitov, A.Walz, 1997 \cite{8}).

In the present paper, we prove that flexible polyhedra do exist in the Min\-kow\-ski
3-space and that each of them preserves the (generalized) volume and
the (total) mean curvature during  a flex.

\section{Existence}

We use standard notations from the theory of the Minkowski spaces freely.
The reader can find them in \cite{13} or \cite{14}.

Let $eu:{\bf R}^n_1\to {\bf R}^n$ be the identity mapping of the Minkowski
$n$-space ${\bf R}^n_1$ to the Euclidean $n$-space ${\bf R}^n$.
Note that a set $D\subset {\bf R}^n_1$
is convex if and only if the set $eu(D)\subset {\bf R}^n$ is convex.
A convex polyhedron is said to be strictly convex if there is no straight
angle among its dihedral angles.

It is crucial for the following Lemma that, according to our definition,
all faces of a polyhedron are triangular.

{\bf Lemma 1.} {\it
Let $P: \Sigma\to {\bf R}^3_1$ be a sphere-type strictly convex polyhedron
and let $Q$ be a disk-type polyhedron which is obtained from
$P$ by removing two adjanced faces. Then $Q$ is flexible.}

{\bf Proof.}
Choose a simplex $\Delta\subset\Sigma$ such that $P(\Delta)\subset Q$.

Let $\Sigma$ have $v$ vertices (0-simplices) $V_1, V_2, \dots , V_v$ and let
$V_{v-2}, V_{v-1}, V_{v}\in\Delta$. Suppose $\Sigma$ has $e$ edges
(1-simplices)  and the edges indexed by $e-2$, $e-1$, and $e$ belong to
$\Delta$.

Consider the family of all polyhedra ${\cal P}:\Sigma\to E^3_{2,1}$ such that
${\cal P}|_\Delta = P|_\Delta$. This family depends on $3v-9$ parameters, $x_1$,
$y_1$, $z_1$, \dots , $x_{v-3}$, $y_{v-3}$, $z_{v-3}$, where $(x_j,
y_j,z_j)={\cal P}(V_j)$ $(j=1,2,\dots, v-3)$ are the coordinates of the $j$th vertex
of $\cal P$.

Suppose vertices ${\cal P}(V_j)$ and ${\cal P}(V_k)$ $(j,k= 1,2,\dots ,v-3)$ are
joined by an edge  and suppose that the edge is indexed by $m$ $(m=1,2,\dots
, e-3)$. Define an auxiliary real-valued function $f_m$ by the formula
$$
f_m(x_1,y_1,z_1,\dots,x_{v-3},y_{v-3},z_{v-3})=
[(x_j-x_k)^2+(y_j-y_k)^2-(z_j-z_k)^2]/2.
$$
The components of the vector-valued function
$f=(f_1,f_2,\dots,f_{e-3})$ represent the half of the squared edge lengths of all
edges of $\cal P$ which do not belong to ${\cal P}(\Delta)$.

Euler's formula yields $3v-e=6$ and, thus, $e-3=3(v-3)$. This means that the Jacobian
matrix of $f$ is a square matrix. Obviously, its $m$th row is as follows:
$$
\matrix {(0 & \dots & 0 & x_j-x_k & y_j-y_k & z_k-z_j & 0 & \dots}
$$
$$
\matrix {\dots & 0 & x_k-x_j & y_k-y_j & z_j-z_k & 0 & \dots & 0)}.
$$

The determinant of this Jacobian matrix does not vanish at the point
corresponding to $P: \Sigma\to {\bf R}^3_1$. To prove this statement, we have to
repeat the above constructions for the convex polyhedron $eu(P)\subset
{\bf R}^3$.
As a result, we obtain a vector-valued function $g=(g_1,g_2,\dots,g_{e-3})$ whose
components
$$
g_m(x_1,y_1,z_1,\dots,x_{v-3},y_{v-3},z_{v-3})=
[(x_j-x_k)^2+(y_j-y_k)^2+(z_j-z_k)^2]/2
$$
$(m=1,2,\dots, e-3)$ represent the half of the squared (Euclidean) lengths of those edges of $eu({\cal P})$
which do not belong to $eu({\cal P})(\Delta)$. The $m$th row of the Jacobian matrix
of $g$ is as follows:
$$
\matrix{(0 & \dots & 0 & x_j-x_k & y_j-y_k & z_j-z_k & 0 & \dots}
$$
$$
\matrix{\dots & 0 & x_k-x_j & y_k-y_j & z_k-z_j & 0 & \dots & 0)}.
$$
This means that the determinant of the Jacobian matrix of $f$ is the product of
$\pm 1$ by the determinant of the Jacobian matrix of $g$, while the latter is
known to be nonzero at the point corresponding to the convex polyhedron $eu(P)$
of the Euclidean 3-space. This statement was obtained in \cite{9} by direct
calculations and is known to be equivalent to the first-order rigidity of a
strictly convex polyhedron in the Euclidean 3-space \cite{1}.

Since the determinant of the Jacobian matrix of $f$ does not vanish at the point
corresponding to $P: \Sigma\to {\bf R}^3_1$, it follows that $f$ maps a
neighborhood $U$ of that point homeomorphically onto its image $f(U)$.
For $t$ close enough to zero, there exists a point $f^t=(f_1^t, f_2^t,\dots,
f_{e-3}^t)\in f(U)$  such that

$\alpha$) $f_m^t$ equals the sum of $t$ and the value of $f_m$ at the point
corresponding to $P$, if $m$ corresponds to the edge of $\Sigma$ shared by those
two simplices which should be removed from $P$ to obtain $Q$,

and

$\beta$) $f_m^t$ equals the value of $f_m$ at the point corresponding to $P$, otherwise.

Since $f$ is a local homeomorphism, it follows that there exists a polyhedron
$P_t^*$ such that $f$ maps the point of $U$ corresponding to $P_t^*$ into $f^t$. This
means that, arbitrarily close to $P$, there exists a polyhedron $P_t^*$ which
is not congruent to $P$ and shares with $P$ the lengths of all edges but the edge
of $\Sigma$ shared by those two simplices which are removed from $P$ to
obtain $Q$.
Removing from $P_t^*$ those two simplices, we obtain a polyhedron $Q_t^*$ which is arbitrary close to $Q$, is
isometric to $Q$, and is not congruent to $Q$.
In \cite{11} it is shown that existence of such a family of polyhedra $Q_t^*$ implies
existence of a family of polyhedra $Q_t$ which is analytic with respect to the
parameter $t$ and satisfies the following conditions: i) $Q_0=Q$; ii) the length
of each edge of $Q_t$ is independent of $t$; iii) there exist two vertices of
$Q_t$ such that the (spatial) distance between them is nonconstant in $t$.
Hence, the family $Q_t$ gives a nontrivial flex of $Q$.
\framebox{\phantom{o}}

Let $\alpha$ be either a time- or a space-like straigt line.
For each $x\in {\bf R}^3_1$, draw a (Minkowski) perpendicular $\alpha^\perp$
to $\alpha$ through $x$ and denote the intersection point  of
$\alpha$ and $\alpha^\perp$ by $y$.
Let $f_\alpha (x)$ be such a point that $f_\alpha\in \alpha^\perp$ and $|f_\alpha
(x)-y|=|x-y|=|f_\alpha (x)-x|/2$. Such a mapping $f_\alpha :{\bf R}^3_1 \to {\bf
R}^3_1$ is said to be a reflection in $\alpha$.

Note that the reflection in $\alpha$ is an isometry of the Minkowski space (it suffices
to verify this statement for $\alpha$ being either the first or third coordinate
axis while in these cases $f_\alpha$ is determined by the matrices
$$
\pmatrix {-1 &  0 &  0\cr
           0 & -1 &  0\cr
           0 &  0 &  1\cr}
\qquad{\rm or}\qquad
\pmatrix { 1 &  0 &  0\cr
           0 & -1 &  0\cr
           0 &  0 & -1\cr}
$$
respectively and, obviously, is a (Minkowski) isometry.)

We need the following statement whose Euclidean version may be found in \cite{10}:

{\bf Lemma 2.} {\it
Let $q_1q_2q_3q_4$ be a spatial quadrilateral in ${\bf R}^3_1$ with the
following properties (see Fig. 1):

1) the opposite sides of $q_1q_2q_3q_4$ are pairwise equal, i.e.,
$|q_1-q_2|=|q_3-q_4|$ and $|q_2-q_3|=|q_4-q_1|$;

2) none of the diagonals $q_1q_3$ and $q_2q_4$ lies on the nullcone;

3) the middle point $q_5$ of the diagonal $q_1q_3$ does not coinside with the
middle point $q_6$ of the diagonal $q_2q_4$;

4) the straigh line $\alpha$ which passes through $q_5$ and $q_6$ does not lie
on the nullcone.

Then the reflection in $\alpha$ maps the quadrilateral $q_1q_2q_3q_4$
onto itself.}

\begin{center}
\quad
\epsfbox{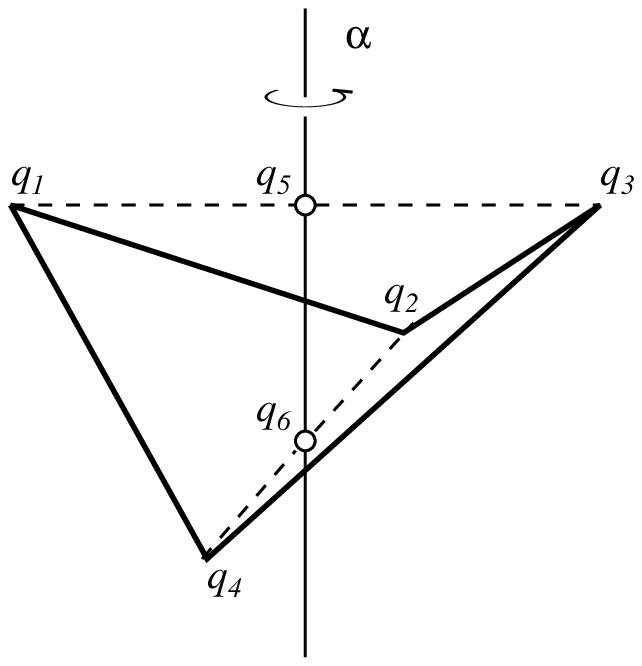}
\quad
\end{center}

\begin{center}
Fig. 1
\end{center}

{\bf Proof.}
It follows from the hypothesis that the triangles $q_1q_2q_3$ and
$q_1q_3q_4$ have pairwise equal sides.
Evidently, the parallelogram equation $|x+y|^2+|x-y|^2=2|x|^2+2|y|^2$
holds true for arbitrary vectors $x,y\in{\bf R}^3_1$.
Hence, the length of the median $(x+y)/2$ can be expressed in terms
of the lengths of the edges $x$, $y$, and $x+y$ of the triangle.
Thus the (corresponding) medians of the triangles $q_1q_2q_3$ and
$q_1q_3q_4$  have equal lengths: $|q_2-q_5|=|q_4-q_5|$. Hence, the
triangle $q_2q_4q_5$ is isosceles and the segment $q_5q_6$ is not only its median
but also its altitude. (The latter can be verified directly: if
$x=(x_1,x_2,x_3)$ and $y=(y_1,y_2,y_3)$ represent equal sides of an isosceles
triangle in ${\bf R}^3_1$, $|x|=|y|$, then $(x+y)/2$ and $x-y$ represent a median
and the third side of the triangle, and by direct calculations we have
$(x+y,x-y)=|x|^2-|y|^2=0$ regardless of the signature of the scalar product.)
Thus we conclude that $\alpha$ is orthogonal to the straight line $q_2q_4$ and, therefore,
the reflection in $\alpha$ interchanges the points $q_2$ and $q_4$.

Similarly, the triangles $q_1q_2q_4$ and $q_2q_3q_4$ have pairwise equal sides,
their medians $q_1q_6$ and $q_3q_6$ have equal lengths, $\alpha$ is orthogonal
to the straight line $q_1q_3$, and the reflection in $\alpha$ interchanges
the points $q_1$ and $q_3$. As a result, we conclude that this reflection maps
the quadrilateral $q_1q_2q_3q_4$ onto itself.
\framebox{\phantom{o}}

{\bf Theorem 1.} {\it
There exist flexible sphere-type polyhedra in the Minkowski 3-space.}

{\bf Proof.}
Let $Q$ be a disk-type polyhedron in ${\bf R}^3_1$ such that

A) $Q$ is obtained from some strictly convex sphere-type polyhedron $P$
by removing some two adjacent (triangular) faces;

B) the boundary $\partial Q$ of $Q$ is a quadrilateral that satisfies the
conditions of Lemma 2;

C) $Q$ is not symmetric with respect to the line of symmetry of $\partial Q$.

From Lemma 1 it follows that $Q$ is flexible. Let $Q_t$ $(0\leq t\leq 1)$
represents a nontrivial flex of $Q$.
Using Lemma 2, we conclude that, for every $t$, the
boundary $\partial Q_t$ of $Q_t$ is symmetric with respect to the reflection in
some straight line. Denote that reflection by $R_t$.
Then $Q_t\cup R_t(Q_t)$ $(0\leq t\leq 1)$ is a nontrivial flex of the
sphere-type polyhedron $Q\cup R_0(Q)$.
\framebox{\phantom{o}}

{\bf Remark 1.}
The examples constructed in Theorem 1 are non-embedded and,
in a sense, they resemble those given by R.~Bricard in the
Euclidean case \cite{Bricard}.
The problem remains open whether there are embedded flexible polyhedra in
the Minkowski 3-space.

\section{Volume}

We define the volume of a domain $D$ in the Minkowski $n$-space ${\bf R}^n_1$
$(n\geq 2)$ as the (Euclidean) volume of its image $eu(D)$ in the Euclidean
$n$-space ${\bf R}^n$, i.e., we put by definition
$$
{\rm vol}\, D=\int_{eu(D)} dx_1dx_2\cdots dx_n.
$$

This definition goes back at least as far as \cite{14} and is natural, since it
introduces a function which is additive, $n$th order homogeneous, and invariant
under isometries of ${\bf R}^n_1$.

{\bf Lemma 3.} {\it
Let $p_0, p_1, \dots , p_n$ be $n+1$ points in the Minkowski space
${\bf R}^n_1$, let $d^2_{jk}=|p_j-p_k|^2$ $(j\neq k=0,1,\dots ,n)$ be the
squared pairwise distances, and let $[p_0,p_1, \dots ,p_n]$ stand for the
simplex with vertices $p_0, p_1, \dots , p_n$. Then
$$
{\rm vol}^2[p_0,p_1,\dots,p_n]=
{(-1)^n\over 2^n(n!)^2}\det
\pmatrix{
0     & 1         & 1        & 1        & \cdot & 1          \cr
1     & 0         & d_{01}^2 & d_{02}^2 & \cdot & d_{0n}^2 \cr
1     & d_{10}^2  & 0        & d_{12}^2 & \cdot & d_{1n}^2 \cr
1     & d_{20}^2  & d_{21}^2 & 0        & \cdot & d_{2n}^2 \cr
\cdot & \cdot     & \cdot    & \cdot    & \cdot & \cdot      \cr
1     & d_{n0}^2  & d_{n1}^2 & d_{n2}^2 & \cdot & 0          \cr}.
\eqno(1)
$$
}

Note that the right-hand side of (1) differs  only in the `minus' sign from the
well-known Cayley--Menger formula representing the volume of a Euclidean simplex
via its edge-lengths \cite{5}.

{\bf Proof.}
Letting the coordinates of $p_j$ be $( p_j^{(1)},p_j^{(2)},\dots,p_j^{(n)})$
$(j=0, 1, \dots , n)$ and using our definition of the volume, from a formula of
elementary analytic geometry we obtain
$$
{\rm vol}\, [p_0,p_1,\dots,p_n]=
{1\over n!}\det
\pmatrix{
p_0^{(1)}  & p_0^{(2)}  & \cdot  & \cdot & p_0^{(n)} & 1      \cr
p_1^{(1)}  & p_1^{(2)}  & \cdot  & \cdot & p_1^{(n)} & 1      \cr
 \cdot     & \cdot      & \cdot  & \cdot & \cdot     & \cdot  \cr
p_n^{(1)}  & p_n^{(2)}  & \cdot  & \cdot & p_n^{(n)} & 1      \cr}.
$$

The determinant is unaltered in value by boarding it with $(n+2)$th row and
column, with `intersecting' element 1, and the remaining elements zero.
Multiplying this boarded determinant by the transpose of the determinant
obtained from it by interchanging the first two rows and last two columns and by
multiplying its $n$th column by $(-1)$, we have
$$
{\rm vol}^2 [p_0,p_1,\dots,p_n]=
{1\over (n!)^2}\det
\pmatrix{
(p_0, p_0)   & (p_0,p_1)    & \cdot  & \cdot & (p_0,p_n)    & 1      \cr
(p_1, p_0)   & (p_1,p_1)    & \cdot  & \cdot & (p_1,p_n)    & 1      \cr
\cdot\ \cdot & \cdot\ \cdot & \cdot  & \cdot & \cdot\ \cdot & \cdot  \cr
(p_n,p_0)    & (p_n,p_1)    & \cdot  & \cdot & (p_n,p_n)    & 1      \cr
1            & 1            & \cdot  & \cdot & 1            & 0      \cr},
\eqno(2)
$$
where
$(p_j,p_k)=p_j^{(1)}p_k^{(1)} + p_j^{(2)}p_k^{(2)} + \cdots +
p_j^{(n-1)}p_k^{(n-1)} - p_j^{(n)}p_k^{(n)}$ $(j,k=0,1,\dots , n)$ stands for
the scalar product in ${\bf R}^n_1$.

If we substitute $(p_j,p_k)=[(p_j,p_j)+(p_k,p_k)-d^2_{jk}]/2$
$(j,k=0,1,\dots ,n)$ in the determinant in (2), subtract from the $j$th row the
product of the last row by $(p_{j-1},p_{j-1})/2$ $(j=1,2,\dots ,n)$, we obtain
(1) after easy reductions. \framebox{\phantom{o}}

Let $P$ be a compact orientable polyhedron in ${\bf R}^n_1$ without boundary, and let $p_*\in
{\bf R}^n_1$.
We define the generalized volume of $P$ as
$$
{\rm Vol}\, P=\sum\varepsilon(p_*,\Delta){\rm vol}\, [p_*, \Delta].
$$
Here the sum is taken over all positively oriented simplices $\Delta$
of $P$,  $[p_*,\Delta]$ stands for the $n$-dimensional simplex that is
the convex hull of $p_*\cup\Delta$ with orientation  generated by that of
$\Delta$, and $\varepsilon(p_*,\Delta)$ is equal to $+1$ (respectively, $-1$)
if the orientation of $[p_*, \Delta]$ agrees with (respectively, is opposite
to) the orientation of the whole space ${\bf R}^n_1$.

Note that the value of the generalized volume is independent of the choice of
$p_*$ (because it is known to be independent for a Euclidean space)
and, for $P$ an embedded polyhedron, ${\rm Vol}\, P$ equals the
volume ${\rm vol}\, D$ of the bounded domain $D$ having $P$ as its boundary.

{\bf Theorem 2.} {\it
The generalized volume of a flexible polyhedron in the Min\-kow\-ski 3-space remains
constant during a flex.}

{\bf Proof.}
The proof repeats literally any of the proofs of the similar theorem in the
Euclidean 3-space given in \cite{8} or [15--17]. The only change is needed:
an extra `minus' sign should be added in all formulas derived from the
Cayley--Menger determinant (1). \framebox{\phantom{o}}

\section{Oriented angle}

Recall from elementary geometry that a number $\varphi_0$ ($0\leq
\varphi_0\leq\pi$)
is called the angle between two nonzero vectors $x,y\in{\bf R}^2$ if
$\varphi_0$
equals the doubled area of the smallest sector that is the intercept of
the unit sphere by the vectors (i.e. directed segments) $x/|x|$ and $y/|y|$.

However, as is known, it is more convenient to treat the
(oriented) angle as a multivalued function $\varphi = \varphi_1 + 2\pi
n$ ($n\in {\bf Z}$) which satisfies the relation
$$
\cos \varphi = {(x,y)\over |x|\cdot |y|},
\eqno(3)
$$
where either $\varphi_1=\varphi_0$ or $\varphi_1=\pi- \varphi_0$.

We use a similar approach to the notion of angle in the Minkowski plane ${\bf R}^2_1$.
The nullcone divides  ${\bf R}^2_1$ into four sectors (connected
components).
Denote them by $S_1,\dots , S_4$ (we assume that $(1,0)\in S_1$,
$(0,1)\in S_2$, $(-1,0)\in S_3$, and $(0,-1)\in S_1$).
According to the classical definition, the angle between two non-null
nonzero vectors  $x,y\in {\bf R}^2_1$  is defined only if $x,y\in S_j$
for some $j=1,\dots ,4$.
In that case, the absolute value of the angle $\theta_0$ equals the
bounded area of the sector that is the intercept of the unit (or
imaginary-unit, depending on $j$) circle by the vectors $x/|x|$ and
$y/|y|$.
We also accept the following `sign convention':
given a positively oriented ordered pair of vectors $x,y\in S_j$
($j=1,\dots,4$), the sign of the angle $\theta_0$ between $x$ and $y$
is positive if $j=1$ or $j=3$ and is negative if $j=2$ or $j=4$. This
definition may be found, for example, in \cite{11} where, in particular,
it is shown that the following relation between the angle and
the scalar product holds:
$$
\cosh \theta_0 = {(x,y)\over |x|\cdot |y|},
\eqno(4)
$$
and, for every sector $S_j$ ($j=1,\dots, 4$) and every real number
$\alpha$,
there exist vectors $x,y\in S_j$ such that the (oriented) angle
between $x$ and $y$ equals $\alpha$.

Following \cite{10}, we extend the above classical definition and treat the angle
as a multivalued function satisfying (4).
Namely, suppose $x,y\in S_1$.
Then the (oriented) angle $\theta_0$ between $x$ and $y$ is already defined.
Rotate $y$ in the positive direction.
The angle $\theta_0$ increases and tends to plus infinity as $y$
tends to the nullcone.

As soon as $y$ intersects the nullcone for the first time
(and, thus, is located in $S_2$), the fraction in (4)
becomes negative, and we put by definition the (oriented) angle between
$x\in S_1$ and $y\in S_2$ equal to the (uniquely determined) complex number
$\theta = \theta_0-i\pi/2$ ($\theta_0\in{\bf R}$) for which the equality
$$
\cosh \theta = {(x,y)\over |x|\cdot |y|}
\eqno(5)
$$
holds true. Note that, while $y$ is rotated in the positive direction within the
limits of $S_2$, the real part $\theta_0$ of the angle $\theta =
\theta_0-i\pi/2$ decreases from minus infinity to plus infinity.

After $y$ intersects the nullcone for the second time, the fraction in (5)
becomes positive again, and we put by definition the (oriented) angle between $x\in
S_1$ and $y\in S_3$ equal to the complex number $\theta = \theta_0-i\pi$
($\theta_0\in \bf R$) for which (5) holds true.
Note that, generally speaking,  there exist two vectors $y, \tilde{y}\in S_3$ of the
prescribed length for which the fraction in (5) takes the same value. Similarly,
there exist two real numbers $\theta_0$ and $\tilde{\theta}_0$ such that $\cosh
(\theta_0-i\pi)=\cosh (\tilde{\theta}_0-i\pi)$. We adopt the following `sign
convention': if the ordered pair $y$ and $\tilde y$ is positively oriented and
if $\theta >\theta_0$ then we put by definition the (oriented) angle
between $x$ and $y$ equal to $\theta$ and the angle between $x$ and $\tilde y$
equal to $\tilde\theta$. Note that, according to the above definition, the real
part $\theta_0$ of the angle $\theta =\theta_0-i\pi$ increases from minus
infinity to plus infinity while $y$ is rotated in the positive direction within
the limits of $S_3$.

As soon as $y$ intersects the nullcone for the third time, the fraction in (5)
becomes negative again, and we put by definition the (oriented) angle
between $x\in S_1$ and $y\in S_4$ equal to the (uniquely determined) complex
number $\theta = \theta_0-i3\pi/2$ ($\theta_0\in \bf R$) for which (5) holds
true. Clearly, while $y$ is rotated in the positive direction within the limits
of $S_4$, the real part $\theta_0$ of the angle $\theta = \theta_0-i3\pi/2$
between $x\in S_1$ and $y\in S_4$ decreases from plus infinity to minus
infinity.

We assume by definition that, each time when, being rotated in the positive
direction, $y$ intersects the nullcone, the imaginary part of the angle receives
the additional summand $-i\pi/2$.

If $x\in S_1$ and $y$ is any non-null nonzero vector, we put by definition
the (oriented) angle between the ordered pair of vectors $y,x$ equal to minus the
(oriented) angle between the ordered pair $x,y$.

We can illustrate the above definition of the angle between a vector $x\in S_1$
and some (non-null nonzero) vector $y\in{\bf R}^2_1$ by Fig. 2.

\begin{center}
\quad
\epsfbox{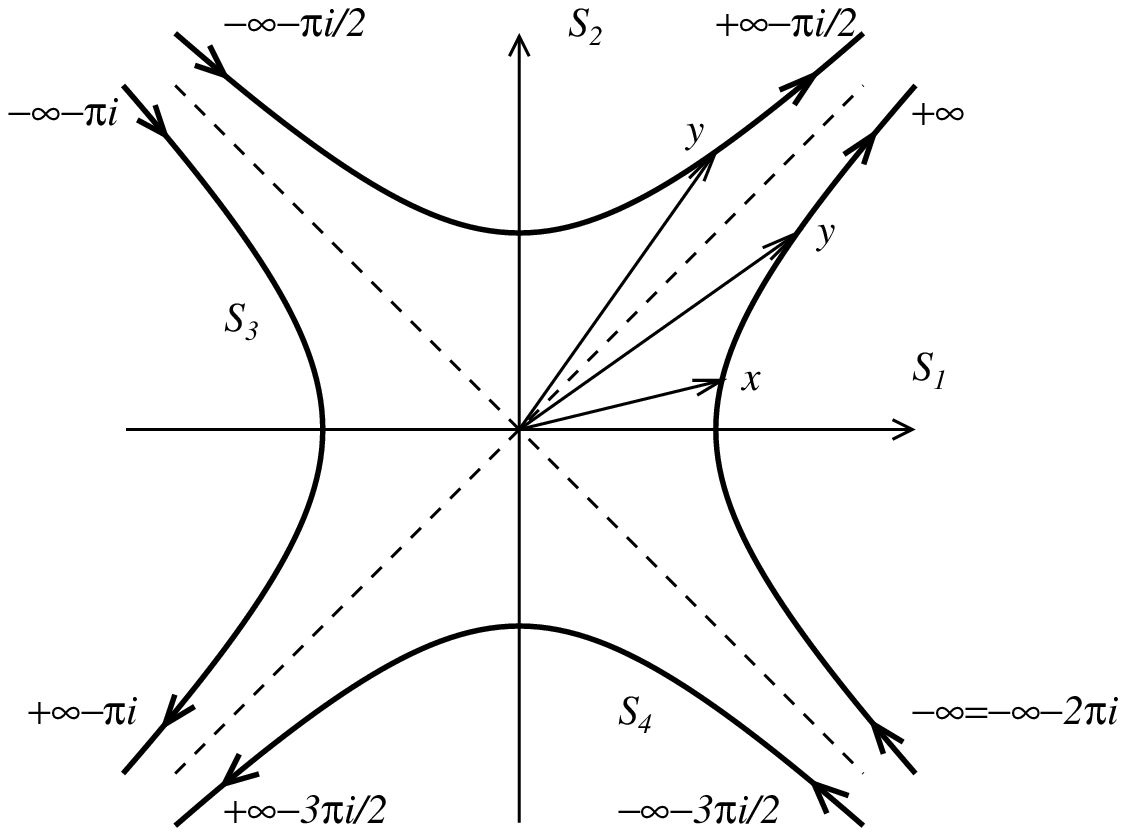}
\quad
\end{center}

\begin{center}
Fig. 2
\end{center}

By definition the (oriented) angle between two non-null nonzero vectors
$x,y\in{\bf R}^2_1$ is the sum of the (oriented) angle between the
vectors $x$ and $e=(1,0)$ and the (oriented) angle between the vectors $e$ and
$y$; denoted by $\angle xy$.

The next four lemmas are needed for the sequel.

{\bf Lemma 4.} {\it
The angle is an  additive function (i.e., if $x,y,z$ are three non-null nonzero vectors and
$\angle xy=\theta_1+ik_1\pi/2$ and $\angle yz=\theta_2+ik_2\pi/2$ then there exists $n\in\bf Z$ such that
$\angle xz=(\theta_1+\theta_2)+i(k_1+k_2)\pi/2 + 2\pi in$).}

{\bf Proof.}
By definition, we have $\angle xy=\angle xe +\angle ey$ and $\angle
yz=\angle ye +\angle ez$. Taking it into account that $\angle ey=-\angle ye$ and
summing the above equations, we obtain $\angle xy+\angle yz =\angle xe +\angle ez
=\angle xz$. \framebox{\phantom{o}}

A nonzero vector $y\in {\bf R}^2_1$ is said to be a right normal vector  to a non-null
nonzero vector $x\in {\bf R}^2_1$ if $y$ is perpendicular to $x$ (i.e., if
$(x,y)=0$) and the ordered pair of vectors $x,y$ is positively oriented.

{\bf Lemma 5.} {\it
Suppose $x\in {\bf R}^2_1$ is a non-null nonzero vector and $y$ is a right
normal vector to $x$. Then $\angle xy = -i\pi/2$.}

{\bf Proof.}
Consider the case $x=(x_1,x_2)\in S_1$. Without loss of generality,
we may assume that $x_1^2-x_2^2=1$ and $y=(x_2,x_1)$. Hence, $y\in S_2$ and thus
$\angle ey =\theta -i\pi/2$ for some real $\theta$. Then $\cosh\angle ey=
-i\sinh\theta =-ix_2$ and $\sinh\angle ey=-i\cosh\theta=\pm\sqrt{-x_2^2-1} =\pm
ix_1$. As soon as $\cosh\theta >0$ and $x_1>0$, we should choose the `minus' sign in
the last formula.

Similarly, $x\in S_1$ and thus $\angle ex=\varphi$ for some real number
$\varphi$. Then $\cosh\angle ex= \cosh\varphi = x_1$ and $\sinh\angle ex=
\sinh\varphi = \pm\sqrt{x_1^2-1}=\pm x_2$. According to our `sign convention',
$\varphi$ should increase as soon as $x_2$ increases; thus, we should choose the
`plus' sign in the last formula.

Finally we have $\cosh\angle xy = \cosh (\angle xe +\angle ey)= \cosh\angle xe
\cdot \cosh\angle ey + \sinh\angle xe \cdot \sinh\angle ey = x_1 (-ix_2) +
(-x_2) \cdot (-ix_1) = 0$ and $\cosh\angle xy =\cosh (\theta - i\pi/2
-\varphi)=-i\sinh (\theta -\varphi)$. Hence, $\theta
=\varphi$ and $\angle xy = -i\pi/2$. This proves the lemma for the case $x\in
S_1$.

Consider the case $x=(x_1,x_2)\in S_2$. Without loss of generality,
we may assume that $x_1^2-x_2^2=-1$ and $y=-(x_2,x_1)$. Hence, $y\in S_3$ and thus
$\angle ey =\theta -i\pi$ for some real $\theta$. Then $\cosh\angle ey=
-i\cosh\theta =-x_2$ and $\sinh\angle ey=-\sinh\theta=\pm\sqrt{x_2^2-1} =\pm
x_1$. According to our `sign convention', $\theta$ should tend to plus infinity
as $x_1\to +\infty$. Thus, we should choose the `minus' sign in the
last formula.

From the above calculations for the case $x\in S_1$ we immediately obtain $\angle ex
=\varphi -i\pi/2$ ($\varphi\in\bf R$), $\cosh\angle ex = -ix_1$, and $\sinh\angle ex = -ix_2$.

Finally, we have $\cosh\angle xy = \cosh (\angle xe +\angle ey)= (-ix_1)\cdot (-x_2)
+(ix_2)\cdot (-x_1) = 0$ and $\cosh\angle xy =\cosh (\theta - i\pi
-\varphi + i\pi/2)=-i\sinh (\theta -\varphi)$. Hence, $\theta
=\varphi$ and $\angle xy = -i\pi/2$. This proves the lemma for the case $x\in
S_2$.

Consider the case $x=(x_1,x_2)\in S_3$. Without loss of generality,
we may assume that $x_1^2-x_2^2=1$ and $y=(x_2,x_1)$. Hence, $y\in S_4$ and thus
$\angle ey =\theta -i3\pi/2$ for some real $\theta$. Then $\cosh\angle ey=
i\sinh\theta =-ix_2$ and $\sinh\angle ey=i\cosh\theta=\pm i\sqrt{x_2^2+1} =\pm i
x_1$. As soon as $\cosh\theta >0$ and $x_1<0$, we should choose the `minus' sign in
the last formula.

From the above calculations for the case $x\in S_2$ we immediately obtain $\angle ex
=\varphi -i\pi$ ($\varphi\in\bf R$), $\cosh\angle ex = x_1$, and $\sinh\angle ex = x_2$.

Finally, we have $\cosh\angle xy = \cosh (\angle xe +\angle ey) = x_1 (-ix_2)
+(-x_2)\cdot (-ix_1) = 0$ and $\cosh\angle xy =\cosh (\theta -
i3\pi/2 -\varphi + i\pi)=-i\sinh (\theta -\varphi)$. Hence, $\theta
=\varphi$ and $\angle xy = -i\pi/2$. This treats the case $x\in S_3$.

Consider the case $x=(x_1,x_2)\in S_4$. Without loss of generality,
we may assume that $x_1^2-x_2^2=1$ and $y=-(x_2,x_1)$. Hence, $y\in S_1$ and thus
$\angle ey =\theta\in \bf R$. Then $\cosh\angle ey=-x_2$ and $\sinh\angle ey=\pm
i\sqrt{x_2^2-1} =\pm x_1$. According to our `sign convention', $\theta$ should
tend to plus infinity as $x_1\to -\infty$. Thus, we should choose the
`minus' sign in the last formula.

From the above calculations for the case $x\in S_3$ we immediately obtain $\angle ex
=\varphi -i3\pi/2$ ($\varphi\in\bf R$), $\cosh\angle ex = -ix_1$, and $\sinh\angle ex = -ix_2$.

Finally, we have $\cosh\angle xy = \cosh (\angle xe +\angle ey)= (-ix_1)\cdot (-x_2)
+(ix_2)\cdot (-x_1) = 0$ and $\cosh\angle xy =\cosh (\theta
-\varphi + i3\pi/2)=-i\sinh (\theta -\varphi)$. Hence, $\theta
=\varphi$ and $\angle xy = -i\pi/2$. This proves the lemma for the case $x\in
S_4$.  \framebox{\phantom{o}}

Let $x,y\in {\bf R}^2_1$ be non-null nonzero vectors and either $|y|=1$ or
$|y|=i$. A real number $t$ is said to be the orthogonal projection of $x$ to
the oriented line spanned by $y$ if the vector $x-ty$ is orthogonal to $y$,
i.e., we put by definition $t=(x,y)$ if $y$ is spacelike, and $t=-(x,y)$ if $y$
is timelike.

{\bf Lemma 6.} {\it
Let $x,y\in {\bf R}^2_1$ be non-null nonzero vectors and either $|y|=1$ or
$|y|=i$. Let $t$ be the orthogonal projection of $x$ to the oriented line
spanned by $y$. Then $t=|x|\cosh\angle xy$ for spacelike $y$ and
$t=-i|x|\cosh\angle xy$ for timelike $y$.}

{\bf Proof.}
By direct calculations
$$
\cosh\angle xy = {(x,y) \over {|x|\cdot |y|}}=
{(x-(x-ty),y) \over {|x|\cdot |y|}}=
t{|y|^2 \over {|x|\cdot |y|}}. \qquad \rm{\framebox{\phantom{o}}}
$$

{\bf Lemma 7.} {\it
Let an ordered pair $a,b\in{\bf R}^2_1$ be positively oriented and such
that $|a|=1$, $|b|=i$. Then, for an arbitrary non-null vector $x\in {\bf R}^2_1$, we
have $x=a|x|\cosh\angle ax +b|x|\sinh\angle ax$ (in other words, if $a,b$ is a
coordinate basis then $x=(|x|\cosh\angle ax ,|x|\sinh\angle ax)$).}

{\bf Proof.}
Let $t$ be the orthogonal projection of $x$ to the oriented line
spanned by $a$ and let $s$ be the orthogonal projection of $x$ to the oriented
line spanned by $b$. Then $x=ta+sb$. Lemma 6 implies $t=|x|\cosh\angle ax$.
By Lemmas 4 and 5, $\angle bx=\angle ba +\angle ax =-\angle ab +\angle ax =
i\pi/2 +\angle ax$. Now, by Lemma 6, $s=-i|x|\cosh\angle bx= -i |x|\cosh
(i\pi/2 +\angle ax)=|x|\sinh\angle ax$.
\framebox{\phantom{o}}

Note that the value of the oriented angle is invariant under orientation preserving
isometries of ${\bf R}^2_1$.

\section{Mean curvature}

Let $P$ be an $(n-1)$-dimensional orientable polyhedron in the Euclidean $n$-space
${\bf R}^n$ ($n\geq 2$) and let $F$ be the set of $(n-1)$-dimensional faces of
$P$. Let $m_j$ ($j\in F)$ be the inward pointing normal to the
$(n-1)$-dimensional faces of $P$, with length equal to the $(n-1)$-volume of the
corresponding $(n-1)$-face. According to the well-known observation of
Minkowski, we have $\sum_{j\in F} m_j=0$ (see,  for example, \cite{1}--\cite{2}).

Lemma 8 below represents a similar result for the Minkowski plane. Note that, in
that case, we should use inward normals for spacelike edges and outward normals for
timelike edges (or vice versa). We use the following notation. For
$w=(w_1,w_2)\in {\bf R}^2_1$, we put by definition $\varepsilon (w)={\rm sgn}\, |w|^2
={\rm sgn}\, (w^2_1-w^2_2)$ and denote the modulus of the complex number $|w|$ by
$\| w\|$. A vector $w$ is said to be unit if $\| w\| =1$.

{\bf Lemma 8.} {\it
Let $g_j\in{\bf R}^2_1$ ($j=1,2,\dots ,k$) be non-null vectors such that
$\sum_{j=1}^k g_j=0$. For each $j=1,2,\dots ,k$, let $n_j$ denote the right
normal unit vector to $g_j$. Then $\sum_{j=1}^k \varepsilon (g_j)\| g_j\| n_j=0$.}

{\bf Proof.}
Let $g_j=(x_j,y_j)$ ($j=1,2,\dots ,k$). Put $u_j=(y_j,x_j)$ by definition.
By assumption, $\sum_{j=1}^k g_j=(\sum_{j=1}^k x_j, \sum_{j=1}^k y_j ) = (0,0)$.
We prove that $\varepsilon (g_j)\| g_j\| n_j=u_j$.

If $g_j$ is spacelike then $u_j$ is a right normal to $g_j$, since
$(u_j,g_j)=0$ and
$$
\det\pmatrix {x_j & y_j \cr
          y_j & x_j \cr}=|g_j|^2>0.
$$
Thus, the vectors $\varepsilon\| g_j\| n_j$ and $u_j$ lie on the same ray
and coincide (because they have the same length).

If $g_j$ is timelike then $-u_j$ is a right normal to $g_j$, since
$(-u_j,g_j)=0$ and
$$
\det\pmatrix {x_j & y_j \cr
          -y_j & -x_j \cr}=-|g_j|^2>0.
$$
Thus, the vectors $\varepsilon (g_j)\| g_j\| n_j$ and $u_j$ lie on the same ray and
coincide.

Finally, $\sum_{j=1}^k \varepsilon (g_j)\| g_j\| n_j =\sum_{j=1}^k u_j= (\sum_{j=1}^k y_j,
\sum_{j=1}^k x_j )=0.$
\framebox{\phantom{o}}

Introduce the notion of the nonoriented angle between two
non-null nonzero vectors $x,y\in{\bf R}^n_1$ ($n\geq 2$). Let $\Pi$ be an oriented
2-dimensional plane passing through $x$ and $y$. Then the following holds:

$\bullet$ \quad
If $\Pi$ is determined uniquely and is spacelike then we order the vectors $x,y$
in such a way that the oriented angle (with respect to $\Pi$) between them is
equal to $\varphi_0+2\pi k$ ($0\leq \varphi_0\leq\pi$, $k\in {\bf Z}$). In this
case, $\varphi_0$ is said to be the nonoriented angle between $x$ and $y$.

$\bullet$ \quad
If $\Pi$ is determined uniquely and is timelike then we order the vectors $x,y$
in such a way that the real part $\theta_0$ of the oriented angle (with respect
to $\Pi$) between $x,y$ is positive. In this case, $\theta_0$ is said to be the
nonoriented angle between $x$ and $y$.

$\bullet$ \quad
If $\Pi$ is determined uniquely and carries a degenerate metric then the
nonoriented angle between $x$ and $y$ is not defined.

$\bullet$ \quad
If $x$ and $y$ lie on a straight line $l$ then we put by definition the
nonoriented angle between $x$ and $y$ equal to 0 for $l$ being spacelike and
$(x,y)>0$, equal to $\pi$ for $l$ being spacelike and $(x,y)<0$, and equal to 0 for
$l$ being timelike.

Let $P\subset{\bf R}^3_1$ be a closed orientable polyhedron such that each edge
of $P$ is non-null and each face carries a nondegenerate metric. Let $E$ and $F$
stand for the sets of edges and faces of $P$. Suppose $f_1,f_2\in
F$ have $g\in E$ as a common edge. Denote by $m_j$ the outward pointing unit
normal to $f_j$ ($j=1,2$). Denote by $\theta (g)$ the nonoriented angle between
$m_1$ and $m_2$. The number
$$
M(P)= {1 \over 2} \sum_{g\in E} \theta (g) \varepsilon (g) \| g\|
\eqno(6)
$$
is called the (total) mean curvature of $P$.

Note that if the above definition is applied to a polyhedron in the Euclidean
3-space, we arrive at the usual definition of the (total) mean curvature of a
polyhedron (cf. \cite{2}).

Let $P(t): \Sigma \to {\bf R}^3_1$ ($0\leq t\leq 1$) be a smooth family of
closed orientable polyhedra such that, for each $t$, each edge of $P(t)$ is
non-null and each face carries a nondegenerate metric. Let $E$ and $F$
stand for the sets of edges and faces of $P(t)$. Suppose $f_1(t), f_2(t)\in
F$ have $g(t)\in E$ as a common edge. Denote by $n_j(t)$ ($j=1,2$) the unit
vector which lies on the plane spanned by $f_j(t)$, is perpendicular to the edge
$g(t)$, and is pointed inward the face $f_j(t)$. Denote by $m_j(t)$ ($j=1,2$) the
outward pointing unit normal to the face $f_j$ of the oriented polyhedron $P(t)$.

{\bf Lemma 9.} {\it
Under the above notation,
$$
{d\theta (g(t))\over dt} = \biggl( {dm_1\over dt}, n_1(t)\biggr)+ \biggl(
{dm_2\over dt}, n_2(t)\biggr)
\eqno(7)
$$
for all $t$ such that $m_1(t)\neq\pm m_2(t)$.}

{\bf Proof.}
Let $\Pi (t)$ be an oriented 2-dimensional plane which is
orthogonal to the edge $g(t)$. Since $g(t)$ does not lie on the nullcone, $\Pi
(t)$ carries a nondegenerate metric. Consider two cases separately.

Case I: $g(t)$ is timelike. Then $\Pi (t)$ carries a Euclidean metric. Let
$e_1, e_2$ be an orthonormal positively oriented frame in $\Pi (t)$.
Denote by $\varphi_j(t)$ ($j=1,2$) the oriented angle between the vectors $e_1$
and $m_j(t)$. If need be, interchange the indicies $j=1$ and $j=2$ in such a
way that the equation $\theta (g(t))= \varphi_1(t) - \varphi_2(t)$ holds true
for $t$ and all real numbers sufficiently close to $t$ (see Fig. 3). Then

\begin{center}
\quad
\epsfbox{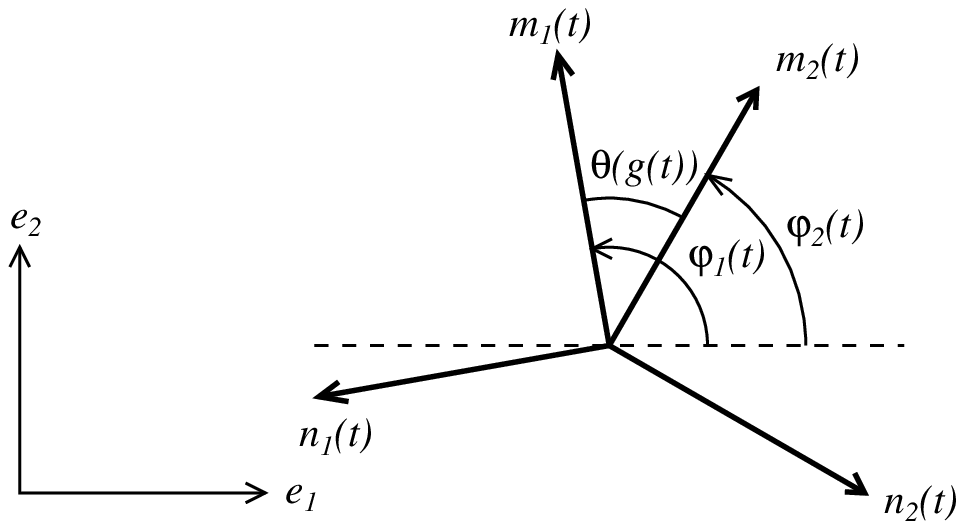}
\quad
\end{center}

\begin{center}
Fig. 3
\end{center}

$$
m_1(t)=\cos\varphi_1(t) e_1 +\sin\varphi_1(t) e_2,
$$
$$
m_2(t)=\cos\varphi_2(t) e_1 +\sin\varphi_2(t) e_2,
$$
$$
n_1(t)=-\sin\varphi_1(t) e_1 +\cos\varphi_1(t) e_2,
$$
$$
n_2(t)=\sin\varphi_2(t) e_1 -\cos\varphi_2(t) e_2,
$$
$$
{dm_1\over dt}=[-\sin\varphi_1(t) e_1 +\cos\varphi_1(t) e_2]{d\varphi_1\over dt},
$$
$$
{dm_2\over dt}=[-\sin\varphi_2(t) e_1 +\cos\varphi_2(t) e_2]{d\varphi_2\over dt},
$$
$$
\biggl( {dm_1\over dt}, n_1(t)\biggr)+ \biggl({dm_2\over dt}, n_2(t)\biggr)=
[\sin^2\varphi_1(t) +\cos^2\varphi_1(t)]{d\varphi_1\over dt}
$$
$$
+[-\sin^2\varphi_2(t) -\cos^2\varphi_2(t)]{d\varphi_2\over dt}=
{d\varphi_1\over dt} - {d\varphi_2\over dt} = {d\theta (g(t))\over dt}.
$$
This proves (7) for the case under consideration.

Case II: $g(t)$ is spacelike. Then $\Pi (t)$ carries a Minkowski metric. Let
$e_1, e_2$ be the standard positevely oriented unit frame in $\Pi (t)$ (in
particular, $e_1$ is spacelike). Denote by $\varphi_j(t)$ ($j=1,2$) the oriented
angle between the vectors $e_1$ and $m_j(t)$. If need be, interchange the
indicies $j=1$ and $j=2$ in such a way that $\theta (g(t))$ equals the real
part of the difference $\varphi_1(t) - \varphi_2(t)$ for $t$ and all real
numbers sufficiently close to $t$ (see Fig. 4). Then, according to Lemma 7,

\begin{center}
\quad
\epsfbox{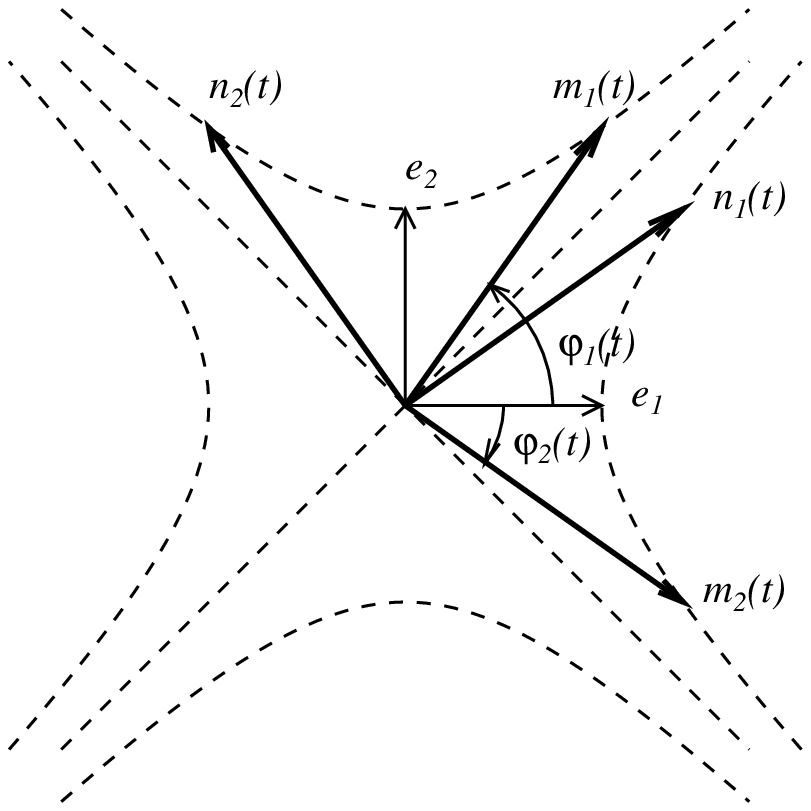}
\quad
\end{center}

\begin{center}
Fig. 4
\end{center}

$$
m_1(t)=|m_1|\cosh\varphi_1(t) e_1 + |m_1|\sinh\varphi_1(t) e_2,
$$
$$
m_2(t)=|m_2|\cosh\varphi_2(t) e_1 + |m_2|\sinh\varphi_2(t) e_2.
$$
Using the usual formulas of hyperbolic trigonometry, $\cosh (\varphi \pm
i\pi/2)= \pm i\sinh\varphi$, $\sinh (\varphi \pm i\pi/2)= \pm i\cosh\varphi$, and
taking into account $|n_j|= -i\varepsilon (n_j) |m_j|$ ($j=1,2$), we obtain
$$
n_1(t)= |n_1|\cosh(\varphi_1 + i\pi/2)e_1 + |n_1|\sinh(\varphi_1 + i\pi/2)e_2
$$
$$
=\varepsilon (n_1) |m_1|\sinh\varphi_1 e_1 + \varepsilon (n_1)
|m_1|\cosh\varphi_1 e_2,
$$
$$
n_2(t)= |n_2|\cosh(\varphi_2 - i\pi/2)e_1 + |n_2|\sinh(\varphi_2 - i\pi/2)e_2
$$
$$
=-\varepsilon (n_2) |m_2|\sinh\varphi_2 e_1 - \varepsilon (n_2)
|m_2|\cosh\varphi_2 e_2,
$$
$$
{dm_1\over dt}= [|m_1|\sinh\varphi_1 e_1 + |m_1|\cosh\varphi_1
e_2]{d\varphi_1\over dt},
$$
$$
{dm_2\over dt}= [|m_2|\sinh\varphi_2 e_1 + |m_2|\cosh\varphi_2
e_2]{d\varphi_2\over dt},
$$
$$
\biggl( {dm_1\over dt}, n_1(t)\biggr)+ \biggl({dm_2\over dt}, n_2(t)\biggr)=
\varepsilon(n_1) |m_1|^2 [\sinh^2\varphi_1-\cosh^2\varphi_1]{d\varphi_1 \over
dt}
$$
$$
- \varepsilon(n_2) |m_2|^2 [\sinh^2\varphi_2-\cosh^2\varphi_2]{d\varphi_2
\over dt} = {d\varphi_1 \over dt} - {d\varphi_2 \over dt} = {d\theta (g(t))
\over dt}.
$$
To derive the last formula, we used the relations $\varepsilon (n_j) |m_j|^2 = -1$
($j=1,2$). This completes the proof of (7) in case II.
\framebox{\phantom{o}}

{\bf Theorem 3.} {\it
Let $P_t$ ($0\leq t\leq 1$) be a closed orientable flexible polyhedron such
that, for each $t$, each edge of $P_t$ is non-null and each face carries a
nondegenerate metric. Then the total mean curvature of $P_t$ is constant during
a flex.}

{\bf Proof.} Let $E$ and $F$ stand for the sets of edges and faces of $P_t$.
Given a face $f(t)\in F$, denote by $m(t)$ the outward
pointing normal unit vector to $f(t)$. Given a face $f(t)\in F$ and an edge
$g(t)\subset f(t)$, denote by $n(t)$ the unit vector which lies on the plane
spanned by $f(t)$, is perpendicular to the edge $g(t)$, and is pointed inward
to the face $f(t)$.

Taking it into account that the edge lengths are independent in $t$, from (6) we obtain
for all $t$ such that there is no edge $g(t)$ of $P_t$ with $\theta(g(t))=0$
or $\theta(g(t))=\pi$:
$$
{d\over dt} M(P_t)={1\over 2}\sum_{g\in E}
{d\theta (g(t)) \over dt} \varepsilon (g(t)) \| g(t)\|
$$
now we use Lemma 9 and sum, first, over all edges of a given face and, then,
over all faces:
$$
=\sum_{f\in F}\sum_{g\in E; g\subset f} \biggl({dm\over dt}, n\varepsilon (g) \|
g\|\biggr)
=\sum_{f\in F} \biggl({dm\over dt}, \sum_{g\in E; g\subset f}
n\varepsilon (g) \| g\|\biggr)
$$
now we apply Lemma 8:
$$
=\sum_{f\in F} \biggl({dm\over dt},0\biggr)=0. \eqno(8)
$$

Let $\tau\in [0,1]$. If there is an edge $g(\tau)$ of $P_\tau$ with
$\theta(g(\tau))=0$ or $\theta(g(\tau))=\pi$, replace one of the faces sharing
$g(\tau)$ by the lateral surface of a tetrahedron whose base is the given face
and whose altitude is small enough in such a way that the resulting polyhedron
$Q_\tau$ has no edges with dihedral angles which are equal to 0 or
$\pi$. The polyhedron $Q_\tau$ is flexible and there is an open
interval $(a,b)$, containing $\tau$, such that, for all $t\in (a,b)$, there is
no edge $g(t)$ of $Q_t$ such that $\theta(g(t))=0$ or $\theta(g(t))=\pi$.
From (8) it follows that $M(Q_t)$ is constant in $t$ on $(a,b)$. Obviously,
$M(P_t)$ is a linear combination of the total mean curvatures of $Q_t$ and the
above tetrahedra each of which moves as a rigid body. Hence, $M(P_t)$ is constant
on $(a,b)$ and thus on $[0,1]$. \framebox{\phantom{o}}

{\bf Remark 2.}
In [18, 19], certain notions of the nonoriented angle between two non-null
nonzero vectors $x,y\in{\bf R}^n_1$ are introduced which are different from
that used in this article.

{\bf Acknowledgement}
The author was partially supported by INTAS-RFBR grant no. IR--97--1778.
He also wishes to thank Institut de Math\'ematiques de Jussieu
(Universit\'e Paris 7), where this paper was written, for
its hospitality and Harold Rosenberg, Idjad Sabitov, and Rabah
Souam for fruitful discussions.

\end{document}